\documentclass[12pt]{amsart}
\usepackage{amsmath,txfonts}
\usepackage{amssymb}
\usepackage{amsxtra}
\usepackage{amsthm, color}
\usepackage{txfonts}
\usepackage{graphicx}
\usepackage{times}
\usepackage{tikz}
\usepackage{hyperref}
\hypersetup{pdfborder={0 0 0.5}}
\usepackage{stmaryrd}
\usepackage[T3,T1]{fontenc}
\DeclareSymbolFont{tipa}{T3}{cmr}{m}{n}
\DeclareMathAccent{\invbreve}{\mathalpha}{tipa}{16}

\usepackage{pgfplots}

\allowdisplaybreaks

\makeatletter
\@namedef{subjclassname@1991}{Subject}
\@namedef{subjclassname@2000}{Subject}
\@namedef{subjclassname@2010}{Subject}
\@namedef{subjclassname@2020}{{\rm 2020} Mathematics Subject Classification}
\makeatother

\def\rr{{\mathbb R}}
\def\rn{{{\rr}^n}}

\def\XXint#1#2#3{{\setbox0=\text{$#1{#2#3}{\int}$ }
		\vcenter{\text{$#2#3$ }}\kern-.6\wd0}}

\newtheorem{theorem}{Theorem}[section]
\newtheorem{definition}[theorem]{Definition}

\newtheorem{proposition}[theorem]{Proposition}
\newtheorem{lemma}[theorem]{Lemma}

\theoremstyle{definition}
\newtheorem{remark}[theorem]{Remark}

\newcommand{\lt}[1]{[ {#1}] \lower.3ex\text{$_{t}$}}

\begin{document}
	\title[Overdetermined affine elliptic problems and characterizations of ellipsoids]
	{Overdetermined affine elliptic problems and characterizations of ellipsoids}

	\author{Jin Dai}

        \address{School of Mathematical Sciences, Chongqing University of Technology, Chongqing, 400054, China}
	            \email{daijin@cqut.edu.cn}

	\author{Tuo Wang}

        \address{School of Mathematics and Statistics, Shaanxi Normal University, Xi'an, 710119, China}
	            \email{wang.tuo@snnu.edu.cn}



	\thanks{JD was supported by the Scientific Research Foundation of Chongqing University of Technology \#2023ZDZ037;
TW was supported by the National Natural Science Foundation of China (Grant No. 12471055), the Recruitment Program for Young Professionals of China, the Youth Innovation Team of Shaanxi Universities, the Shaanxi Fundamental Science Research Project for Mathematics and Physics (Grant No. 22JSZ012), and the Fundamental Research Funds for the Central Universities (Grant Nos. GK202603001, GK202307001, and GK202202007). }

	\date{\today}

	\subjclass[2020]{35N25, 52A40}

	\keywords{Affine $p$-Laplacian, Overdetermined problem, Alexandrov theorem}

	\begin{abstract} In this work, we mainly establish an affine analogue of Serrin's symmetry theorem \cite{Serrin}. And we also establish an affine Alexandrov theorem and an inequality stronger than the Petty projection inequality in convex geometry.
				\end{abstract}

	\maketitle

\tableofcontents
	
	\arraycolsep=1pt
	\numberwithin{equation}{section}

\section{Introduction} \label{s1}
The study of overdetermined problems began with the celebrated paper by Serrin \cite{Serrin}. In particular, he proved the radial symmetry of the solution to the following overdetermined elliptic problem. Let $\Omega$ be a bounded domain in $\rn$ with smooth boundary, let $\nu_\Omega$ be the outward unit normal to $\partial\Omega$, and let $c$ be a constant. If $u$ is a smooth solution to
\begin{equation*}
\left\{
  \begin{array}{ll}
    \Delta u=-1,~& \mathrm{in}~\Omega; \\
    u=0,~& \mathrm{on}~\partial\Omega;\\
    \frac{\partial u}{\partial \nu_\Omega}=c,~& \mathrm{on}~\partial\Omega,
  \end{array}
\right.
\end{equation*}
then $\Omega$ must be a ball and, up to a translation, $u(x)=\frac{R^2-{|x|}^2}{2n}$ for some $R>0$. Here, $|x|$ denotes the standard Euclidean norm of $x\in\rn$.

Motivated by the seminal works of Zhang \cite{Zhang99} and Lutwak--Yang--Zhang \cite{LYZ02}, Haddad--Jim\'enez--Montenegro \cite{Haddad} introduced the affine $p$-Laplace operator in divergence form (the case $p=2$ was previously defined by Schindler--Tintarev \cite{ST18}) by
$$\Delta_p^\mathcal{A}u=\mathrm{div}(H_u^{p-1}(\nabla u)\nabla H_u(\nabla u)),$$
for $p\geq 1$ and a function $u\in W^{1,p}_0(\Omega)\backslash\{0\}$.
Here, $H_u$ is defined by
\begin{equation}\label{44}
  H_u^p(v)=c_{n,p}^{-n}\mathcal{E}_p^{p+n}(u)\int_{S^{n-1}}\bigg(\int_\rn|\nabla u\cdot\xi|^pdx\bigg)^{-\frac{n+p}{p}}|v\cdot\xi|^pd\xi,\ \ \ \forall v\in\rn,
\end{equation}
where $S^{n-1}$ denotes the unit sphere, $d\xi$ is its standard surface measure, $\cdot$ denotes the Euclidean inner product, $c_{n,p}$ is the normalizing constant (see \eqref{38}), and the affine energy is defined by
\begin{equation}\label{52}
  \mathcal{E}_p(u)=c_{n,p}\bigg(\int_{S^{n-1}}\bigg(\int_\rn|\nabla u\cdot\xi|^pdx\bigg)^{-\frac{n}{p}}d\xi\bigg)^{-\frac{1}{n}},
\end{equation}
and $W^{1,p}_0(\Omega)=\big\{u\in W^{1,p}(\rn):u=0~\text{a.e. on}~\rn\backslash\Omega\big\}$ is the subspace of the Sobolev space $W^{1,p}(\rn)$.

Notice that $H_u$ is a norm of class $C^2$ on $\rn\backslash\{0\}$ (see Lemmas \ref{11} and \ref{46}), so the affine $p$-Laplace operator is exactly the Finsler $p$-Laplace operator, that is,
$$\Delta_p^\mathcal{A}u=\Delta_p^{H_u}u.$$
For a sufficiently regular norm $H:\rn\rightarrow[0,+\infty)$, the Finsler $p$-Laplace operator (see \cite{BFK}) is defined by
$$\Delta_p^Hu=\mathrm{div}(H^{p-1}(\nabla u)\nabla H(\nabla u))$$
with $p\geq1$. The overdetermined anisotropic problem
\begin{equation}\label{30}
\left\{
  \begin{array}{ll}
    \Delta_p^H u=-1,~& \mathrm{in}~\Omega; \\
    u=0,~& \mathrm{on}~\partial\Omega;\\
    H(\nabla u)=c,~& \mathrm{on}~\partial\Omega
  \end{array}
\right.
\end{equation}
was solved for $p>1$ (see \cite{BC18,LLZ}), and $\Omega$ is characterized, up to translations and dilations, by the Wulff shape of $H$, that is,
$$\overline{\Omega}=\bigcap_{y\in\rn}\Big\{x\in\rn:x\cdot y\leq H(y)\Big\}.$$
This leads us to the problem of finding $u$ and a constant $c>0$ satisfying the following non-local overdetermined system for $p>1$:
\begin{equation}\label{6}
\left\{
  \begin{array}{ll}
    \Delta_p^\mathcal{A} u=-1,~& \mathrm{in}~\Omega; \\
    u=0,~& \mathrm{on}~\partial\Omega;\\
    H_u(\nabla u)=c,~& \mathrm{on}~\partial\Omega.
  \end{array}
\right.
\end{equation}
The additional boundary condition in \eqref{6} is understood pointwise. For domains with $C^{2,\alpha}$ boundary, this interpretation is justified by the $C^{1,\alpha}(\overline{\Omega})$ regularity established in Section \ref{s3}, which ensures that $\nabla u$ extends continuously to $\partial\Omega$.

For the non-local Dirichlet problem
\begin{equation}\label{1}
\left\{
  \begin{array}{ll}
    \Delta_p^\mathcal{A} u=-1,~& \mathrm{in}~\Omega; \\
    u=0,~& \mathrm{on}~\partial\Omega,
  \end{array}
\right.
\end{equation}
we adopt the following notion of weak solution:
\begin{definition}
Let $p>1$. We say that $u$ is a weak solution of \eqref{1} if $u\in W^{1,p}_0(\Omega)\backslash\{0\}$ and
$$\int_\Omega H_u^{p-1}(\nabla u)\nabla H_u(\nabla u)\cdot \nabla\phi dx=\int_\Omega\phi dx\ \ \ \forall\phi\in C^1_0(\Omega).$$
\end{definition}
A solution of the non-local overdetermined problem \eqref{6} is understood to be a weak solution in the above sense.

Following \cite{Haddad}, we consider the affine energy functional
$$E^\mathcal{A}_p(u)=\frac{\mathcal{E}_p^p(u)}{\|u\|_1^p}$$
for $p>1$ and $u\in W^{1,p}_0(\Omega)\backslash\{0\}$.
We show that minimizers of the variational problem
\begin{equation}\label{2}
  \frac{1}{\tau_p^\mathcal{A}(\Omega)}=\inf_{u\in W^{1,p}_0(\Omega)\backslash\{0\}}E^\mathcal{A}_p(u)
\end{equation}
satisfy the associated Euler--Lagrange equation and, after rescaling, yield solutions of \eqref{1}:
\begin{lemma}\label{10}
Let $\Omega$ be a bounded open set and let $p>1$. Then the infimum of the variational problem \eqref{2} is attained.
Consequently, there exists a weak solution to the non-local Dirichlet problem \eqref{1}.
\end{lemma}
For a fixed solution $u$, the norm $H_u$ is fixed, so equation \eqref{1} can be regarded as the anisotropic Dirichlet problem
\begin{equation}\label{3}
\left\{
  \begin{array}{ll}
    \Delta_p^{H_u} u=-1,~& \mathrm{in}~\Omega; \\
    u=0,~& \mathrm{on}~\partial\Omega.
  \end{array}
\right.
\end{equation}
Thus, regularity results for solutions to \eqref{3} yield corresponding regularity results for solutions of the non-local Dirichlet problem \eqref{1}.

Throughout, $C^{2,\alpha}(\rn\backslash\{0\})$ denotes the space of functions whose second derivatives are locally $\alpha$-H\"older continuous on $\rn\backslash\{0\}$. To solve the overdetermined problem, we require additional regularity. For the anisotropic $p$-Laplacian, we consider the class of norms defined by
$$\mathcal{I}_p=\Big\{H\in C^{2,\alpha}(\rn\backslash\{0\})\ \text{for some}\ \alpha\in(0,1):\frac{H^p}{p}\in C^2_+(\rn\backslash\{0\})\Big\},$$
where $C^2_+(\rn\backslash\{0\})$ denotes the class of twice continuously differentiable functions whose Hessian matrices are strictly positive definite at every point.
We now recall Bianchini--Ciraolo's result for the overdetermined anisotropic problem.
\begin{theorem}\label{29}
(\cite{BC18}) Let $\Omega$ be a bounded domain with $C^{2,\alpha}$ boundary and let $H\in \mathcal{I}_p$ with $p>1$. There is a solution $u$ to \eqref{30} for some constant $c>0$ if and only if, up to translations and dilations, $\overline{\Omega}$ is the Wulff shape of $H$ and
there exist $z\in\rn$ and $r>0$ such that $\Omega=\{x\in\rn:H_0(x)<r\}+z$ and
$$u(x)=\frac{(p-1)(r^{\frac{p}{p-1}}-H_0(x-z)^{\frac{p}{p-1}})}{pn^{\frac{1}{p-1}}},$$
for each $x\in\Omega$. For each prescribed $c>0$, the solution is unique when it exists. Here $H_0$ is the dual norm of $H$.
\end{theorem}
In particular, taking $H=h_{B^n}$, the support function of the $n$-dimensional unit ball, gives the solution of the overdetermined problem for the classical $p$-Laplacian:
\begin{theorem}\label{35}
(\cite{BC18,FK}) Let $\Omega$ be a bounded domain with $C^{2,\alpha}$ boundary and let $p>1$. Then there exists a solution $u$, for some constant $c>0$, of the overdetermined problem
\begin{equation*}
\left\{
  \begin{array}{ll}
    \Delta_p u=-1,~& \mathrm{in}~\Omega; \\
    u=0,~& \mathrm{on}~\partial\Omega;\\
    |\nabla u|=c,~& \mathrm{on}~\partial\Omega,
  \end{array}
\right.
\end{equation*}
if and only if $\Omega$ is a ball and there exist $z\in\rn$ and $r>0$ such that
$$u(x)=\frac{(p-1)(r^{\frac{p}{p-1}}-|x-z|^{\frac{p}{p-1}})}{pn^{\frac{1}{p-1}}},$$
for each $x\in\Omega$. For each prescribed $c>0$, the solution is unique when it exists.
\end{theorem}

Another important tool is the $L^p$ projection centroid conjecture raised by Lutwak \cite{lutwak} for $p=1$ and by Lutwak--Yang--Zhang \cite{LYZ00} for general $p$. The case $p=1$ was settled by Milman--Shabelman--Yehudayoff \cite{MSY}, and the case $p>1$ was recently settled by Dai--Wang \cite{DW26}:
\begin{theorem}\label{14}
(\cite{DW26,MSY}) Let $K$ be a convex body containing the origin in its interior and let $p\geq 1$. Then there is a constant $c>0$ such that $\Gamma_p \Pi_p^\circ K=cK$ if and only if $K$ is an origin-symmetric ellipsoid.
\end{theorem}

Using Theorems \ref{29}, \ref{35} and \ref{14}, we solve the non-local overdetermined problem \eqref{6} and obtain a characterization of ellipsoids:
\begin{theorem}\label{32}
Let $\Omega$ be a bounded domain with $C^{2,\alpha}$ boundary and let $p>1$. There exist a constant $c>0$ and a solution $u$ to \eqref{6} if and only if $\Omega$ is an ellipsoid and
there exist $z\in\rn$, $r>0$ and $\psi\in \mathrm{SL}(n)$ such that
$$u(x)=\frac{(p-1)(r^{\frac{p}{p-1}}-|\psi x-z|^{\frac{p}{p-1}})}{pn^{\frac{1}{p-1}}},$$
for each $x\in\Omega$. Moreover, for each such solution $u$, $r$ is uniquely determined, and the pair $(\psi,z)$ is unique up to the simultaneous left action
$$(\psi,z)\longmapsto(Q\psi,Qz),\qquad Q\in\mathrm{SO}(n).$$
\end{theorem}

For $p=1$, if $u$ is sufficiently smooth, its affine 1-Laplacian $\Delta^\mathcal{A}_1 u$ may be interpreted in terms of the curvature of its level sets; see the open problem in \cite{Haddad}.
This curvature interpretation naturally leads to Alexandrov-type rigidity problems. We first recall the anisotropic version of Alexandrov's theorem:
\begin{theorem}\label{4}
(\cite{BCS16}) Let $H$ be a norm on $\rn$ of class $C^2(\rn\backslash\{0\})$ such that $H^2$ is uniformly convex, and let $\partial \Omega$ be a compact connected $C^2$ hypersurface without boundary embedded in Euclidean space. If the anisotropic mean curvature $M^H_{\partial \Omega}$ is constant on $\partial \Omega$, then $\overline{\Omega}$ is the Wulff shape of $H$ up to translations and dilations.
\end{theorem}

To characterize ellipsoids, we modify this formulation using the following non-local affine mean curvature:
\begin{definition}
Let $F$ be a set of finite perimeter with $C^2$ boundary in $\rn$, and let $\Gamma_1\Pi_1^\circ F$ be the $L^1$ centroid body of the polar $L^1$ projection body of $F$. Define the non-local affine mean curvature $M_{\partial F}^\mathcal{A}$ of $\partial F$ by
$$M_{\partial F}^\mathcal{A}(x)=\bigg(\frac{\omega_n}{|\Pi_1^\circ F|}\bigg)^{\frac1n}\frac{\mathrm{div}(\nabla h_{\Gamma_1\Pi_1^\circ F}(\nu_F(x)))}{n-1},\ \ \ \forall x\in \partial F,$$
where $h_{\Gamma_1\Pi_1^\circ F}$ is the support function of $\Gamma_1\Pi_1^\circ F$, $\nu_F$ is the outward unit normal to $F$, and $|\Pi_1^\circ F|$ denotes the volume of $\Pi_1^\circ F$.
\end{definition}
\begin{remark}
Let $\Omega$ be a bounded convex domain with $C^2$ boundary containing the origin in its interior, and let $u$ be the zero extension of $g(1/\rho_{\overline{\Omega}})$, where $g\in C^2([0,1])$ is strictly decreasing and $g(1)=0$. For each regular value $t\in(0,g(0))$ of $u|_{\Omega\setminus\{0\}}$, the non-local affine mean curvature on $\partial[u]_t=\{x\in\rn:u(x)=t\}$ is related to the affine 1-Laplacian by (see Lemmas \ref{46} and \ref{34})
$$M^\mathcal{A}_{\partial[u]_t}(x)=-\frac{1}{n-1}\Delta_1^\mathcal{A}u(x),\ \ \ \forall x\in \partial[u]_t.$$
\end{remark}
We then obtain the following affine Alexandrov theorem:
\begin{theorem}\label{5}
Let $\partial \Omega$ be a compact connected $C^2$ hypersurface without boundary embedded in Euclidean space. If the non-local affine mean curvature $M^\mathcal{A}_{\partial \Omega}$ is constant on $\partial \Omega$, then $\Omega$ is an ellipsoid.
\end{theorem}
Finally, we establish an inequality involving the non-local affine mean curvature that is stronger than the Petty projection inequality and characterizes ellipsoids through its equality case.

\begin{theorem}\label{17}
Let $K$ be a convex body with $C^2$ boundary. Then
\begin{equation}\label{26}
  |K|\leq n\omega_n^{2+\frac1n}|\Pi_1^\circ K|^{-\frac1n}\bigg(\int_{\partial K}M_{\partial K}^\mathcal{A}(x)h_{\Gamma_1\Pi_1^\circ K}(\nu_K(x))d\mathcal{H}^{n-1}(x)\bigg)^{-1},
\end{equation}
with equality if and only if $K$ is an ellipsoid.
\end{theorem}
This inequality follows directly from Minkowski's second inequality. It is stronger than Petty's projection inequality \cite{Petty71}, which states that, for each convex body $K$,
\begin{equation}\label{18}
  |K|^{\frac{n-1}{n}}\leq \omega_n|\Pi_1^\circ K|^{-\frac{1}{n}}.
\end{equation}

The remainder of this paper is organized as follows. In Section \ref{s2}, we provide some necessary preliminaries; in Section \ref{s3}, we solve the affine overdetermined problem; and in Section \ref{s4}, we prove the affine Alexandrov theorem and establish related geometric inequalities.
\section{Preliminaries}\label{s2}

In this section, we collect the basic material needed from convex geometry and Sobolev spaces. For more information, see \cite{EG15,Gardner06,Gruber07,Sch14}.

Throughout this paper, we work in the $n$-dimensional Euclidean space $\rn$, $n\geq 2$, with the standard orthonormal basis $\{e_1,e_2,\cdots,e_n\}$. For $p\geq 1$,
the Lebesgue space $L^p(\rn)$ denotes the class of all Lebesgue measurable functions $f$ with the $L^p$ norm
$${\|f\|}_{p} =\Big(\int_{\rn }|f(x)|^pdx\Big)^{\frac{1}{p}}<\infty$$
and the Sobolev space $W^{1,p}(\rn)$ denotes the class of all functions $f\in L^p(\rn)$ whose weak gradient satisfies $|\nabla f|\in L^p(\rn)$.
Further, if $\Omega\subset\rn$ is a domain, we denote by $W^{1,p}_0(\Omega)$ the subspace of $W^{1,p}(\rn)$ consisting of functions vanishing outside $\Omega$:
$$W^{1,p}_0(\Omega)=\big\{u\in W^{1,p}(\rn):u=0~\text{a.e. on}~\rn\backslash\Omega\big\}.$$
Given a function $f:\rn\rightarrow\mathbb{R}$ and $t\in\mathbb{R}$, the super-level set $\{x\in\rn:f(x)>t\}$ is simply denoted by $[f]_t$.
For a linear transformation $\phi$ on $\rn$, we denote its determinant by $\det\phi$ and its transpose by $\phi^*$. If $\phi$ is invertible, we write $\phi^{-*}=(\phi^{-1})^*$.

For a set $A\subset \rn$ and $\lambda\in\mathbb R$, we write
$$\lambda A=\big\{\lambda a: a\in A\big\},$$
and write $-A$ for $(-1)A$. If $\lambda>0$, this is called a dilation of $A$.
The Minkowski sum of sets $A,B\subset \rn$ is denoted by
$$A+B=\big\{a+b: a\in A,b\in B\big\}.$$
In particular, $A+\{z\}$ is written as $A+z$ for $z\in\rn$.
We denote by $\mathcal{H}^{n-1}$ the $(n-1)$-dimensional Hausdorff measure on $\mathbb{R}^n$. For a set $A\subset\rn$, $\overline{A}$, $\mathrm{int}(A)$ and $\partial A$ denote its closure, interior and boundary, respectively.
The symbols $|A|$ and $P(A)$ denote the volume ($n$-dimensional Lebesgue measure) and the perimeter of $A$, respectively. If $P(A)<\infty$, then $A$ is called a set of finite perimeter. If $\partial A$ is a sufficiently smooth hypersurface, then $P(A)=\mathcal{H}^{n-1}(\partial A)$.

For $k\in\mathbb{N}$, the $k$-dimensional unit ball of $\mathbb{R}^k$ is denoted by $B^k$ and the volume $|B^k|$ is denoted by $\omega_k$. For real $t\geq0$, we extend this notation by $\omega_t=\pi^{t/2}/\Gamma(1+t/2)$, where $\Gamma$ is the Gamma function. For $p\geq 1$, set the constant
\begin{equation}\label{38}
  c_{n,p}=(n\omega_n)^\frac{1}{n}\Big(\frac{n\omega_n\omega_{p-1}}{2\omega_{n+p-2}}\Big)^\frac{1}{p}.
\end{equation}

\subsection{Convex bodies and norms}
A convex body $K\subset\rn$ is a convex compact set with nonempty interior. It is uniquely determined by its support function, defined by
$$h_K(x)=\sup\big\{x\cdot y:y\in K\big\},~\forall x\in\rn.$$
Let $K,L$ be convex bodies, let $\lambda>0$, and let $\phi$ be a linear transformation. Then the support function satisfies
$$h_{\lambda K+L}=\lambda h_K+h_L\ \ \ \mathrm{and}\ \ \ h_{\phi K}(x)=h_K(\phi^*x)~\forall x\in\rn.$$
If $x\in\partial K$, one has $h_K(\nu_K(x))=x\cdot\nu_K(x)$.
If a convex body $K$ contains the origin in its interior, define its gauge function $|\cdot|_K$ and radial function $\rho_K(\cdot)$, respectively, by
$$|x|_K=\inf\big\{\lambda\geq 0:x\in \lambda K\big\},~\forall x\in\rn,$$
and
$$\rho_K(x)=\sup\big\{\lambda\geq 0:\lambda x\in K\big\},~\forall x\in\rn\backslash\{0\}.$$
Then, we have $|x|_K=\frac{1}{\rho_K(x)}$ for $x\in\rn\backslash\{0\}$.

If a convex body $K$ contains the origin in its interior, define its polar body $K^\circ$ by
$$K^\circ=\big\{x\in\rn:x\cdot y\leq 1\ \mathrm{for}\ \mathrm{all}\ y\in K\big\}.$$
Then, we have $h_K(x)=|x|_{K^\circ}$ for $x\in\rn$.

A convex body $K$ is said to be origin-symmetric if $K=-K$. For an origin-symmetric convex body $K$, its support function $h_K$ and gauge function $|\cdot|_K$ are both norms. Conversely, if $H:\rn\rightarrow[0,\infty)$ is a norm, then there is an origin-symmetric convex body $K$ such that $H=h_K$ and $H=|\cdot|_{K^\circ}$. If we define the dual norm $H_0$ of $H$ by
$$H_0(x)=\sup_{y\ne0}\frac{x\cdot y}{H(y)}\ \ \ \forall x\in\rn,$$
then $H_0=h_{K^\circ}$ and $H_0=|\cdot|_{K}.$

Let $K$ be a convex body and $\nu_K:\partial K\rightarrow S^{n-1}$ be
the general Gauss map of $K$. For each Borel set $\omega\subset S^{n-1}$, the inverse spherical image $\nu^{-1}_K(\omega)$ of $\omega$ is the set of all boundary points of $K$ which have an outward unit normal vector belonging to $\omega$.

Given a convex body $K$ containing the origin in its interior and $p\geq 1$, we define its $L^p$ surface area measure $S_{p,K}$ on $S^{n-1}$ by
$$S_{p,K}(\omega)=\int_{\nu^{-1}_K(\omega)}h^{1-p}_K(\nu_K(x))d\mathcal{H}^{n-1}(x),$$
for each Borel set $\omega\subset S^{n-1}$. For $p=1$, we write $S_K=S_{1,K}$.

The class of all convex bodies is usually endowed with the topology induced by the Hausdorff distance. It is known that a sequence $\{K_i\}$ of convex bodies converges to a convex body $K$ with respect to the Hausdorff distance, written $K_i\rightarrow K$, if and only if $h_{K_i}\rightarrow h_K$ uniformly on $S^{n-1}$ as $i\rightarrow\infty$.

\subsection{Projection bodies and centroid bodies}
Let $F\subset\rn$ be a set of finite perimeter with Lipschitz boundary $\partial F$, and denote the outward unit normal to $\partial F$ by $\nu_F(x)$ for $\mathcal{H}^{n-1}$-almost every $x\in\partial F$. Define the projection body $\Pi F$ of $F$ by its support function
$$h_{\Pi F}(x)=\frac{1}{2} \int_{\partial F} |\nu_F(y)\cdot x|d\mathcal{H}^{n-1}(y),\quad x\in\rn.$$
For a convex body $K$ and $\xi\in S^{n-1}$, let $K|_{\xi^\perp}$ denote the orthogonal projection of $K$ onto the hyperplane $\big\{x\in\rn:x\cdot\xi=0\big\}$. Then the $(n-1)$-dimensional volume of this projection is given by
\begin{equation*}
  \mathcal{H}^{n-1}(K|_{\xi^\perp})
  =\frac{1}{2}\int_{\partial K}|\xi\cdot\nu_K(x)|\,
  d\mathcal{H}^{n-1}(x).
\end{equation*}
Thus, the support function of the projection body of a convex body equals the volume of its orthogonal projection.

Let $K$ be a convex body containing the origin in its interior. Define the $L^p$ projection body $\Pi_pK$ of $K$ by its support function
$$h_{\Pi_p K}^p(x)=\frac{1}{d_{n,p}} \int_{\partial K} |\nu_K(y)\cdot x|^ph_K^{1-p}(\nu_K(y))d\mathcal{H}^{n-1}(y),\quad x\in\rn,$$
where $d_{n,p}=\frac{2\omega_{n+p-2}}{\omega_{p-1}}$. More generally, for a set $F$ of finite perimeter with Lipschitz boundary, we write $\Pi_1F:=\omega_{n-1}^{-1}\Pi F$; this agrees with the above definition when $F$ is a convex body. The normalization gives
\begin{equation*}
  \Pi_pB^n=B^n.
\end{equation*}
Using the $L^p$ surface area measure, one also has
$$h_{\Pi_p K}^p(x)=\frac{1}{d_{n,p}} \int_{S^{n-1}} |\xi\cdot x|^pdS_{p,K}(\xi),\quad x\in\rn.$$

The following lemma records translation and dilation properties of projection bodies.
\begin{lemma}\label{41}
(\cite{Gardner06,LYZ00}) If $F$ is a set of finite perimeter with Lipschitz boundary, $K$ is a convex body containing the origin in its interior, $c_1>0$ and $z\in\rn$, then
\begin{enumerate}
  \item $\Pi (F+z)=\Pi F;$
  \item for every $p\geq 1$, $\Pi_p(c_1K)=c_1^{(n-p)/p}\Pi_p K.$
\end{enumerate}
\end{lemma}



For $p\geq1$, let $K$ be a convex body containing the origin in its interior. Define the $L^p$ centroid body $\Gamma_pK$ of $K$ by its support function
$$h_{\Gamma_p K}^p(x)=\frac{1}{b_{n,p}|K|} \int_{K} |y\cdot x|^pdy,\quad x\in\rn,$$
where $b_{n,p}=\frac{\omega_{n+p}}{\omega_2\omega_n\omega_{p-1}}$. The constant $b_{n,p}$ is chosen so that
\begin{equation*}
  \Gamma_pB^n=B^n.
\end{equation*}
The classical centroid body $\Gamma K$ of $K$ is defined by
$$h_{\Gamma K}(x)=\frac{1}{|K|} \int_{K} |y\cdot x|dy,\quad x\in\rn,$$
and satisfies $\Gamma K=b_{n,1}\Gamma_1K$.
The following lemma collects some properties of centroid bodies.
\begin{lemma}\label{43}
(\cite{LYZ00}) If $K$ is a convex body containing the origin in its interior, $\phi\in\mathrm{SL}(n)$, $\lambda>0$ and $p\geq 1$, then
$$\Gamma_p \phi K=\phi\Gamma_p K\ \ \ \mathrm{and}\ \ \   \Gamma_p \lambda K=\lambda\Gamma_p K.$$
\end{lemma}
We also need the following regularity property of the $L^p$ centroid body:
\begin{lemma}\label{11}
(\cite{LYZ00}) If $p\geq1$ and $K$ is a convex body containing the origin in its interior, then the centroid body $\Gamma_p K$ is of class $C^2_+$. In particular, $h_{\Gamma_p K}\in C^2(\rn\backslash\{0\})$.
\end{lemma}

Both the $L^p$ centroid body $\Gamma_pK$ and the $L^p$ projection body $\Pi_pK$ are origin-symmetric.
For convenience, we denote the polar bodies of $\Pi_p K$ and $\Gamma_p K$ by $\Pi_p^\circ K$ and $\Gamma_p^\circ K$, respectively.
\subsection{Affine Laplacians}
Let $f\in W^{1,p}(\rn)\backslash\{0\}$ with $p\geq 1$. Define a convex body $L_{p,f}$ of $f$ by its gauge function
\begin{equation}\label{def1}
  |x|^p_{L_{p,f}}=\int_{\rn} |\nabla f(y)\cdot x|^pdy,\quad x\in\rn.
\end{equation}
The body $L_{p,f}$ has the following affine properties.
\begin{lemma}\label{45}
(\cite{Haddad}) Let $f\in W^{1,p}(\rn)\backslash\{0\}$, $\phi\in \mathrm{SL}(n)$ and $g(x)=f(x+z)$ for some $z\in\rn$. Then,
$$L_{p,f\circ\phi}=\phi^{-1}L_{p,f}\ \ \ \mathrm{and}\ \ \ L_{p,g}=L_{p,f}.$$
\end{lemma}
We use the following representation.
\begin{lemma}\label{46}
(\cite{Haddad}) Let $u\in W^{1,p}(\rn)\backslash\{0\}$. Then,
$$H_u=\Big(\frac{\omega_n}{|L_{p,u}|}\Big)^\frac{1}{n}h_{\Gamma_pL_{p,u}}.$$
\end{lemma}
Consequently, for each fixed $u$, the affine $p$-Laplacian $\Delta^\mathcal{A}_p u$ is the anisotropic $p$-Laplacian $\Delta^{H_u}_p u$, and $H_u$ is a constant multiple of the support function of the centroid body $\Gamma_pL_{p,u}$.

By analogy with the $p$-Laplacian, which equals the Fr\'echet derivative of $-\|\nabla u\|^p_p/p$, the affine $p$-Laplacian can be obtained by differentiating the energy
$-\mathcal{E}^p_p(u)/p$ (see \eqref{52}):
\begin{proposition}\label{53}
(\cite{Haddad}) Let $\Omega$ be a bounded domain and $u\in W^{1,p}_0(\Omega)\backslash\{0\}$.
\begin{enumerate}
\item If $p>1$, then we have the G\^ateaux derivative
$$(-\mathcal{E}^p_p(u)/p)'=\Delta_p^\mathcal{A}u.$$
\item If $p=1$ and $\nabla u\neq0$ a.e. in $\Omega$, then we have the G\^ateaux derivative
$$(-\mathcal{E}_1(u))'=\Delta_1^\mathcal{A}u.$$
\end{enumerate}
\end{proposition}

Next, we list some properties of the affine $p$-Laplacian.
\begin{lemma}\label{8}
(\cite{Haddad}) Let $\phi\in \mathrm{SL}(n)$, $\lambda>0$, $p\geq 1$, and $u\in W^{1,p}(\rn)\backslash\{0\}$. Then,
\begin{enumerate}
  \item $\Delta_p^\mathcal{A}(u\circ \phi)=(\Delta_p^\mathcal{A}u)\circ \phi$,
  \item $\Delta_p^\mathcal{A}(\lambda u)=\lambda^{p-1}\Delta_p^\mathcal{A}u$.
\end{enumerate}
\end{lemma}

\begin{proposition}\label{47}
(\cite{Haddad}) Let $p\geq1$. For any nonzero function $u\in W^{1,p}(\rn)$ that is radial about some point $z\in\rn$, namely $u(x)=f(|x-z|)$ for some real-valued function $f$, we have
$$\Delta_p^\mathcal{A}u=\Delta_pu\ \ \ \mathrm{in}\ \rn,$$
and $H_u=h_{B^n}$.
\end{proposition}

\subsection{Mixed volumes and representations of integrals}
We finally collect the mixed-volume formulas and curvature representations needed for the geometric arguments in Section \ref{s4}. Minkowski's theorem leads to the notion of mixed volumes:
\begin{theorem}
(\cite{Gruber07,Sch14}) For convex bodies $K_1,\dots,K_m\subset\rn$, the mixed volumes $V(K_{i_1},\dots,K_{i_n})$, $1\leq i_1,\dots,i_n\leq m$, are symmetric in their indices and satisfy
$$|\sum_{i=1}^m \lambda_iK_i|=\sum_{i_1,\dots,i_n=1}^mV(K_{i_1},\dots,K_{i_n})\lambda_{i_1}\cdots\lambda_{i_n}\quad \mathrm{for}\ \lambda_1,\dots,\lambda_m\geq0.$$
\end{theorem}

For convex bodies $K,L$, we denote $V(K,\dots,K,L)$ and $V(K,\dots,K,L,L)$ by $V_1(K,L)$ and $V_2(K,L)$, respectively, where $K$ appears $n-1$ times in the first expression and $n-2$ times in the second.

For the mixed volume $V_1(K,L)$, we have the representation formula \cite{Gruber07,Sch14}
\begin{equation}\label{48}
  V_1(K,L)=\frac{1}{n}\int_{S^{n-1}}h_LdS_K=\frac{1}{n}\int_{\partial K}h_L\circ\nu_Kd\mathcal{H}^{n-1}.
\end{equation}

To state the representation formula for $V_2(K,L)$, we first recall the associated anisotropic mean curvature. We say that a convex body $L$ is of class $C^2_+$ if $L$ is of class $C^2$ and all principal curvatures of its boundary are positive.
Let $K$ be of class $C^2$, and let $H$ be a norm of class $C^2(\rn\backslash\{0\})$. The anisotropic mean curvature $M^H_{\partial K}$ is defined by
$$M^H_{\partial K}(x)=\frac{\mathrm{div}(\nabla H(\nu_K(x)))}{n-1},\ \ \ \forall x\in \partial K.$$
In particular, if $L$ is of class $C^2_+$, taking $H=h_L$ gives $M^{h_L}_{\partial K}$.
For the mixed volume $V_2(K,L)$, we have the representation formula \cite{BF,Xia2017}
\begin{equation}\label{49}
  V_2(K,L)=\frac{1}{n}\int_{\partial K}M^{h_L}_{\partial K}h_L\circ \nu_Kd\mathcal{H}^{n-1}.
\end{equation}


\section{The Solution of the Affine Overdetermined Problem} \label{s3}
\subsection{The Existence and Regularity of Solutions}
We begin by establishing the existence and regularity of solutions to the non-local Dirichlet problem \eqref{1}.

\begin{lemma}\label{7}
(\cite{Haddad}) Let $\Omega$ be a bounded open set and $p\geq 1$. Then, for any $f\in W_0^{1,p}(\Omega)\backslash\{0\}$,
$$\mathcal{E}_p(f)\geq C_{n,p}(\Omega)\|f\|_p^{\frac{n-1}{n}}\|\nabla f\|^{\frac{1}{n}}_p\geq D_{n,p}(\Omega)\|f\|_1^{\frac{n-1}{n}}\|\nabla f\|^{\frac{1}{n}}_p,$$
where $C_{n,p}(\Omega),D_{n,p}(\Omega)$ are positive constants.
\end{lemma}

\noindent{\bf Proof of Lemma \ref{10}.}
For $\lambda\in\mathbb{R}\backslash\{0\}$ and $u\in W^{1,p}_0(\Omega)\backslash\{0\}$, since $\mathcal{E}_p(\lambda u)=|\lambda|\mathcal{E}_p(u)$ and $\|\lambda u\|_1=|\lambda|\|u\|_1$, we have $E^\mathcal{A}_p(\lambda u)=E^\mathcal{A}_p(u)$. Let $\{u_k\}\subset W^{1,p}_0(\Omega)$ be a sequence such that $\|u_k\|_1=1$ and $\mathcal{E}_p^p(u_k)\rightarrow 1/\tau^\mathcal{A}_p(\Omega)$. Then $\{u_k\}$ is bounded in $W^{1,p}_0(\Omega)$ by Lemma \ref{7}. By compactness, up to a subsequence, $u_k\rightarrow u$ weakly in $W^{1,p}_0(\Omega)$ and $u_k\rightarrow u$ in $L^p(\Omega)$. Clearly, $\|u\|_1=1$. The lower semicontinuity
$$\mathcal{E}_p(u)\leq\liminf_{k\rightarrow\infty}\mathcal{E}_p(u_k)$$
follows from the proof of \cite[Theorem 3]{Haddad} and leads, since $\|u\|_1=1$, to
$$E^\mathcal{A}_p(u)=\frac{\mathcal{E}_p^p(u)}{\|u\|_1^p}=\mathcal{E}_p^p(u)\leq 1/\tau^\mathcal{A}_p(\Omega).$$
Thus $E^\mathcal{A}_p(u)=1/\tau^\mathcal{A}_p(\Omega)$, and hence $u$ is a minimizer of the variational problem \eqref{2}. Since $\mathcal{E}_p(u)=\mathcal{E}_p(|u|)$, we may replace $u$ by $|u|$ and assume that $u\geq0$.

Set
$$\mu=\frac{1}{\tau_p^\mathcal{A}(\Omega)}=\mathcal{E}_p^p(u)>0.$$
The definition of the infimum gives
$$\mathcal{E}_p^p(v)\geq\mu\|v\|_1^p\qquad\text{for every }v\in W_0^{1,p}(\Omega).$$
We now use the following auxiliary-function argument, which avoids differentiating the $L^1$ norm. For every $\phi\in C_0^1(\Omega)$ and all sufficiently small $|t|$, define
$$g(t)=\mathcal{E}_p^p(u+t\phi)-\mu\left(\int_\Omega(u+t\phi)\,dx\right)^p.$$
Since $u\geq0$ and $\int_\Omega u\,dx=\|u\|_1=1$, for sufficiently small $|t|$ we have
$$0<\int_\Omega(u+t\phi)\,dx\leq\int_\Omega|u+t\phi|\,dx=\|u+t\phi\|_1.$$
In particular, $u+t\phi$ is nonzero, so $\mathcal{E}_p^p(u+t\phi)$ is well defined.
Consequently,
$$g(t)\geq\mathcal{E}_p^p(u+t\phi)-\mu\|u+t\phi\|_1^p\geq0,$$
whereas $g(0)=0$. Thus $g$ attains a local minimum at $t=0$. By Proposition \ref{53}, $g$ is differentiable at $0$, and therefore $g'(0)=0$. Computing the derivative gives
$$p\int_\Omega H_u^{p-1}(\nabla u)\nabla H_u(\nabla u)\cdot\nabla\phi\,dx
=p\mu\left(\int_\Omega u\,dx\right)^{p-1}\int_\Omega\phi\,dx.$$
Since $\int_\Omega u\,dx=1$, it follows that
\begin{equation}\label{9}
\int_\Omega H_u^{p-1}(\nabla u)\nabla H_u(\nabla u)\cdot\nabla\phi\,dx
=\mu\int_\Omega\phi\,dx
\qquad\text{for every }\phi\in C_0^1(\Omega).
\end{equation}
Thus $-\Delta_p^\mathcal{A}u=\mu$ in the weak sense. Finally, set $w=\lambda u$ with $\lambda=\mu^{\frac{1}{1-p}}$. Since $H_{\lambda u}=H_u$ for every $\lambda>0$, we have $H_w=H_u$. By Lemma \ref{8},
$$\Delta_p^\mathcal{A}w=\mu^{-1}\Delta_p^\mathcal{A}u=-1.$$
Since $w\in W_0^{1,p}(\Omega)\backslash\{0\}$, multiplying identity \eqref{9} by $\lambda^{p-1}=\mu^{-1}$ shows that $w$ is a weak solution of the non-local Dirichlet problem \eqref{1}.
{\hfill$\Box$}

We record the following regularity result.
\begin{theorem}[\cite{Haddad}]
Let $u$ be a weak solution of the non-local Dirichlet problem \eqref{1}. Then, the following regularity properties hold:
\begin{enumerate}
  \item $u\in C^{1,\alpha}(\Omega)$ for some $0<\alpha<1$,
  \item $u\in C^{1,\alpha}(\overline{\Omega})$ for some $0<\alpha<1$ provided that $\partial \Omega$ is of class $C^{2,\alpha}$.
\end{enumerate}
\end{theorem}

\subsection{Auxiliary Results}

We next collect several auxiliary results needed in the proof of Theorem \ref{32}.

\begin{lemma}\label{31}
(\cite{LXZ11}) If $K$ is a convex body containing the origin in its interior and $g:(0,\infty)\rightarrow\mathbb{R}$ is strictly decreasing, then for $f:\rn\setminus\{0\}\rightarrow \mathbb{R}$ defined by $f(x)=g(1/\rho_K(x))$,
$$|\nabla f(x)|=\frac{|g'(1/\rho_K(x))|}{|z\cdot\nu_K(z)|}$$
holds for almost all $x\in\rn$, where $z={\rho_K(x)x}\in\partial K$ and $\nu_K(z)$ is the outward unit normal. At points where $\nabla f(x)\neq0$,
$$\nu_K(z)=-\frac{\nabla f(x)}{|\nabla f(x)|}=\frac{\nabla(1/\rho_K)(z)}{|\nabla(1/\rho_K)(z)|}.$$
\end{lemma}

\begin{lemma}\label{34}
Let $\Omega$ be a bounded convex domain containing the origin in its interior, and let $p\geq1$. Let $g\in W^{1,p}([0,1])$ be strictly decreasing with $g(1)=0$. Define $f=g(1/\rho_{\overline{\Omega}})$ in $\Omega\setminus\{0\}$, set $f(0)=g(0)$, and set $f=0$ in $\rn\setminus\Omega$. Then, there is a constant $c>0$ independent of $\Omega$ such that
$$L_{p,f}=c\Pi_p^\circ\overline{\Omega}.$$
\end{lemma}
\begin{proof}
The Lipschitz property of $1/\rho_{\overline{\Omega}}$ on $\Omega$, the chain rule, and the co-area formula imply that $f\in W^{1,p}(\Omega)$. Since $g(1)=0$, the zero extension satisfies $f\in W^{1,p}_0(\Omega)$.
By definition \eqref{def1}, the support function is given by
$$h_{L_{p,f}^\circ}(\xi)=\bigg(\int_\Omega |\nabla f\cdot\xi|^pdx\bigg)^\frac{1}{p}.$$
Since $\partial\Omega=\partial\overline{\Omega}$ and $0<1/\rho_{\overline{\Omega}}<1$ in $\Omega\setminus\{0\}$, the computation in Lemma \ref{31} applies with $K=\overline{\Omega}$ on $\Omega$. Setting $t=g(s)$ and using the co-area formula, we calculate that
\begin{eqnarray*}
  \int_\Omega |\nabla f\cdot\xi|^pdx &=& \int_\Omega |\frac{\nabla f}{|\nabla f|}\cdot\xi|^p|\nabla f|^pdx \\
   &=& \int_0^{g(0)}\bigg(\int_{\partial [f]_t}|\frac{\nabla f}{|\nabla f|}\cdot\xi|^p|\nabla f|^{p-1}d\mathcal{H}^{n-1}\bigg)dt \\
   &=& \int_1^{0}\bigg(\int_{s(\partial \Omega)}|\frac{\nabla f}{|\nabla f|}\cdot\xi|^p|\nabla f|^{p-1}d\mathcal{H}^{n-1}\bigg)g'(s)ds\\
   &=& \int_0^{1}|g'(s)|\bigg(\int_{s(\partial \Omega)}|\nu_{\overline{\Omega}}(z)\cdot\xi|^p\Big(\frac{|g'(1/\rho_{\overline{\Omega}}(x))|}{|z\cdot\nu_{\overline{\Omega}}(z)|}\Big)^{p-1}d\mathcal{H}^{n-1}(x)\bigg)ds\\
   &=&\int_0^{1}|g'(s)|\bigg(\int_{s(\partial \Omega)}|\nu_{\overline{\Omega}}(z)\cdot\xi|^p|g'(s)|^{p-1}h_{\overline{\Omega}}^{1-p}(\nu_{\overline{\Omega}}(z))d\mathcal{H}^{n-1}(x)\bigg)ds\\
   &=&\int_0^{1}|g'(s)|^ps^{n-1}ds\int_{\partial \Omega}|\nu_{\overline{\Omega}}(z)\cdot\xi|^ph_{\overline{\Omega}}^{1-p}(\nu_{\overline{\Omega}}(z))d\mathcal{H}^{n-1}(z).
\end{eqnarray*}
Thus,
$$\bigg(\int_\Omega |\nabla f\cdot\xi|^pdx\bigg)^\frac{1}{p}=\bigg(\int_0^{1}|g'(s)|^ps^{n-1}ds\bigg)^\frac{1}{p}\bigg(\int_{\partial \Omega}|\nu_{\overline{\Omega}}(z)\cdot\xi|^ph_{\overline{\Omega}}^{1-p}(\nu_{\overline{\Omega}}(z))d\mathcal{H}^{n-1}(z)\bigg)^\frac{1}{p}.$$
Then, since $g\in W^{1,p}([0,1])$, there is a constant $c_1>0$ such that
$$h_{L_{p,f}^\circ}(\xi)=c_1h_{\Pi_p\overline{\Omega}}(\xi),$$
that is,
$$L_{p,f}=\frac{1}{c_1}\Pi_p^\circ\overline{\Omega}.$$
\end{proof}

\begin{lemma}\label{lem:centroid-regularity}
Let $p>1$, let $K$ be a convex body containing the origin in its interior, and set $H=h_{\Gamma_pK}$. Then
$$H\in C^{2,\alpha_p}(S^{n-1}),\qquad \alpha_p=\min\{p-1,1\}>0.$$
Consequently, $H\in C^{2,\alpha_p}(\rn\backslash\{0\})$.
\end{lemma}
\begin{proof}
By the definition of $\Gamma_pK$, write
$$F(x):=H^p(x)=c_K\int_K|x\cdot y|^pdy,\qquad c_K=(b_{n,p}|K|)^{-1}.$$
For every $q>-1$, slicing $K$ into line segments parallel to $\xi$ gives the uniform estimate
\begin{equation}\label{eq:centroid-slice}
\sup_{\xi\in S^{n-1}}\int_K|\xi\cdot y|^qdy<\infty.
\end{equation}
Indeed, the $(n-1)$-dimensional volumes of the perpendicular sections are uniformly bounded, and the remaining one-dimensional integral is finite because $q+1>0$. This uniform integrability permits differentiation under the integral sign and yields
\begin{align*}
D_iF(x)&=pc_K\int_K|x\cdot y|^{p-1}\operatorname{sgn}(x\cdot y)y_i\,dy,\\
D_{ij}F(x)&=p(p-1)c_K\int_K|x\cdot y|^{p-2}y_i y_j\,dy.
\end{align*}

First suppose that $1<p<2$ and put $\beta=p-2\in(-1,0)$. We claim that
\begin{equation}\label{eq:centroid-holder}
\int_K\big||x\cdot y|^\beta-|x'\cdot y|^\beta\big|dy
\leq C|x-x'|^{\beta+1},\qquad x,x'\in S^{n-1}.
\end{equation}
It suffices to consider sufficiently close $x,x'$; otherwise the estimate follows from \eqref{eq:centroid-slice}. Write
$$x=ce-sv,\qquad x'=ce+sv,\qquad e\perp v,\qquad |x-x'|=2s,$$
where $e,v$ are unit vectors, $c=\cos(\theta/2)$, $s=\sin(\theta/2)$, and $\theta$ is the angle between $x$ and $x'$. Decompose $y=te+uv+z$, where $t,u\in\mathbb R$ and $z\in(\operatorname{span}\{e,v\})^\perp$ denotes the component in the remaining $n-2$ dimensions. For fixed $z$, set $a=ct=c(e\cdot y)$ and $\rho=s|u|=s|v\cdot y|$. Then
$$x\cdot y=ct-su=a-su,\qquad x'\cdot y=ct+su=a+su,$$
and hence
$$\Delta:=\big||x\cdot y|^\beta-|x'\cdot y|^\beta\big|=\big||a-\rho|^\beta-|a+\rho|^\beta\big|,$$
where the last equality follows by interchanging the two terms when $u<0$.
The estimates below are trivial when $\rho=0$, so assume $\rho>0$. On $|a|\leq2\rho$, we have $\Delta\leq|a-\rho|^\beta+|a+\rho|^\beta$. The change of variables $a=\rho r$ therefore gives
$$\begin{aligned}
\int_{|a|\leq2\rho}\Delta\,da
&\leq\rho^{\beta+1}\int_{-2}^{2}\big(|r-1|^\beta+|r+1|^\beta\big)\,dr\\
&=\frac{2(1+3^{\beta+1})}{\beta+1}\rho^{\beta+1}
\leq C\rho^{\beta+1},
\end{aligned}$$
where the last integral is finite because $\beta>-1$. For $a\geq2\rho$, applying the Lagrange mean value theorem to $t\mapsto t^\beta$ on $[a-\rho,a+\rho]$ gives $\tau_a\in(a-\rho,a+\rho)$ such that
$$\Delta=-2\beta\rho\tau_a^{\beta-1}\leq C\rho a^{\beta-1},$$
since $\tau_a>a/2$ and $\beta-1<0$. By symmetry, the same estimate holds with $|a|$ in place of $a$ when $a\leq-2\rho$. Consequently,
$$\int_{|a|\geq2\rho}\Delta\,da
\leq2C\rho\int_{2\rho}^{\infty}a^{\beta-1}\,da
=\frac{2C}{-\beta}\rho(2\rho)^\beta
\leq C\rho^{\beta+1},$$
because $\beta<0$. Combining the two estimates yields
$$\int_{\mathbb R}\Delta\,da=\int_{|a|\leq2\rho}\Delta\,da+\int_{|a|\geq2\rho}\Delta\,da\leq C\rho^{\beta+1}.$$
Let
$$\mathcal D=\{(u,z):te+uv+z\in K\ \text{for some}\ t\in\mathbb R\}$$
be the bounded projected domain. Integration over $\mathcal D$ therefore gives
$$\int_K\Delta\,dy\leq C\int_{\mathcal D}\rho^{\beta+1}\,dudz=C\int_{\mathcal D}(s|u|)^{\beta+1}\,dudz\leq Cs^{\beta+1},$$
which proves \eqref{eq:centroid-holder}. Since $K$ is bounded and $\beta+1=p-1$, it follows that
$$|D_{ij}F(x)-D_{ij}F(x')|\leq C|x-x'|^{p-1},\qquad x,x'\in S^{n-1}.$$
Thus $F\in C^{2,p-1}(S^{n-1})$.

Next, if $p=2$, then $F$ is a quadratic form. If $p>2$, differentiation once more gives
$$D_{ijk}F(x)=p(p-1)(p-2)c_K\int_K|x\cdot y|^{p-3}\operatorname{sgn}(x\cdot y)y_i y_j y_k\,dy.$$
Since $p-3>-1$, estimate \eqref{eq:centroid-slice} and the boundedness of $K$ show that all third derivatives are uniformly bounded on $S^{n-1}$. Hence the second derivatives are Lipschitz on $S^{n-1}$, namely $F\in C^{2,1}(S^{n-1})$.

Finally, because $K$ is full-dimensional, $F(x)>0$ for every $x\neq0$. Thus $F$ has a positive lower bound on $S^{n-1}$, and smoothness of $t\mapsto t^{1/p}$ away from zero implies that $H=F^{1/p}$ belongs to $C^{2,\alpha_p}(S^{n-1})$. The last assertion follows by positive homogeneity of degree one.
\end{proof}

\begin{proposition}\label{33}
Let $p>1$ and let $f\in W^{1,p}(\rn)$ be nonzero. Then $h_{\Gamma_pL_{p,f}}\in \mathcal{I}_p$.
\end{proposition}
\begin{proof}
Set $H=h_{\Gamma_pL_{p,f}}$ and $K=\Gamma_pL_{p,f}$. Since $f$ is nonzero, $L_{p,f}$ is an origin-symmetric convex body with nonempty interior. By Lemma \ref{11}, $K$ is of class $C^2_+$. Moreover, Lemma \ref{lem:centroid-regularity} gives
$$H\in C^{2,\alpha_p}(S^{n-1}),\qquad \alpha_p=\min\{p-1,1\}>0,$$
and hence $H\in C^{2,\alpha_p}(\rn\backslash\{0\})$. When $p\geq2$, local Lipschitz continuity of the second derivatives implies their local $\alpha$-H\"older continuity for every $0<\alpha<1$, so we may replace $C^{2,1}$ by $C^{2,\alpha}$.

It remains to prove that $H^p/p\in C^2_+(\rn\backslash\{0\})$. Fix $\xi\in S^{n-1}$. Since $K$ is of class $C^2_+$, the standard support-function calculation gives
$$A_\xi:=\nabla^2H(\xi)\big|_{\xi^\perp}>0,\qquad \nabla^2H(\xi)\xi=0.$$
Write $\nabla H(\xi)=H(\xi)\xi+z_\xi$, where $z_\xi\in\xi^\perp$. Differentiating $H^p$ yields
$$\nabla^2(H^p)(\xi)=pH(\xi)^{p-1}\nabla^2H(\xi)+p(p-1)H(\xi)^{p-2}\nabla H(\xi)\otimes\nabla H(\xi).$$
Consequently, for every $q=a\xi+w$ with $a\in\mathbb R$ and $w\in\xi^\perp$,
\begin{align*}
q^T\nabla^2(H^p)(\xi)q
&=pH(\xi)^{p-1}w^TA_\xi w\\
&\quad+p(p-1)H(\xi)^{p-2}\big(aH(\xi)+z_\xi\cdot w\big)^2.
\end{align*}
If $w\neq0$, the first term is positive; if $w=0$ and $a\neq0$, the second term is positive because $p>1$ and $H(\xi)>0$. Hence $\nabla^2(H^p)(\xi)>0$ for every $\xi\in S^{n-1}$. Since $\nabla^2(H^p)$ is positively homogeneous of degree $p-2$, this conclusion extends to $\rn\backslash\{0\}$. Therefore $H^p/p\in C^2_+(\rn\backslash\{0\})$, which proves that $H\in\mathcal I_p$.
\end{proof}

\subsection{Proof of the Affine Serrin Theorem}
\noindent{\bf Proof of Theorem \ref{32}.} Let $p>1$. First, we prove sufficiency.
For each ellipsoid $E$, there are $\phi\in \mathrm{SL}(n)$  and a ball $B$ such that $E=\phi B$. For the ball $B$, by Theorem \ref{35}, the overdetermined problem
\begin{equation*}
\left\{
  \begin{array}{ll}
    \Delta_p u=-1,~& \mathrm{in}~\mathrm{int}(B); \\
    u=0,~& \mathrm{on}~\partial B;\\
    |\nabla u|=c,~& \mathrm{on}~\partial B,
  \end{array}
\right.
\end{equation*}
has a unique radial solution $u$, and there exist $z\in\rn$ and $r>0$ such that
\begin{equation}\label{36}
  u(x)=\frac{(p-1)(r^{\frac{p}{p-1}}-|x-z|^{\frac{p}{p-1}})}{pn^{\frac{1}{p-1}}}.
\end{equation}
Since $u$ is radial about $z$, it follows from Proposition \ref{47} that $\Delta_p^\mathcal{A} u=\Delta_p u$ and $H_u(y)=|y|$.
Thus, $u$ is also a solution of the overdetermined affine problem
\begin{equation}\label{37}
\left\{
  \begin{array}{ll}
    \Delta_p^\mathcal{A} u=-1,~& \mathrm{in}~\mathrm{int}(B); \\
    u=0,~& \mathrm{on}~\partial B;\\
    H_u(\nabla u)=c,~& \mathrm{on}~\partial B,
  \end{array}
\right.
\end{equation}
Setting $x=\phi^{-1} y$ and $w(y)=u(x)=u(\phi^{-1}y)$ for each $y\in\rn$, it follows from Lemmas \ref{45}, \ref{43}, \ref{46}, and \ref{8} that $H_w=\Big(\frac{\omega_n}{|L_{p,u}|}\Big)^\frac{1}{n}h_{\phi\Gamma_pL_{p,u}}$
and $\Delta_p^\mathcal{A} u(x)=\Delta_p^\mathcal{A} w(y)$. Then,
$$H_w(\nabla w(y))=\Big(\frac{\omega_n}{|L_{p,u}|}\Big)^\frac{1}{n}h_{\phi\Gamma_pL_{p,u}}(\phi^{-*}\nabla u(x))=H_u(\nabla u(x)).$$
Combining these identities with equation \eqref{37} and $E=\phi B$ shows that $w$ is a solution of the overdetermined affine problem
\begin{equation*}
\left\{
  \begin{array}{ll}
    \Delta_p^\mathcal{A} w=-1,~& \mathrm{in}~\mathrm{int}(E); \\
    w=0,~& \mathrm{on}~\partial E;\\
    H_w(\nabla w)=c,~& \mathrm{on}~\partial E.
  \end{array}
\right.
\end{equation*}
Moreover, formula \eqref{36} gives
$$w(y)=u(\phi^{-1}y)=\frac{(p-1)(r^{\frac{p}{p-1}}-|\phi^{-1}y-z|^{\frac{p}{p-1}})}{pn^{\frac{1}{p-1}}}.$$
It follows that if $\Omega$ is an ellipsoid, then there are $\psi\in \mathrm{SL}(n)$, $z\in\rn$ and $r>0$ such that the function
$$u(x)=\frac{(p-1)(r^{\frac{p}{p-1}}-|\psi x-z|^{\frac{p}{p-1}})}{pn^{\frac{1}{p-1}}},\qquad x\in\Omega,$$
is a solution of problem \eqref{6}. Here $r$ is unique, and $(\psi,z)$ is unique up to the simultaneous left action $(\psi,z)\mapsto(Q\psi,Qz)$ with $Q\in\mathrm{SO}(n)$. Indeed, if two such pairs give the same function, then $A=\psi'\psi^{-1}$ is an isometry with $z'=Az$; since $\det A=1$, one has $A\in\mathrm{SO}(n)$.

Next, we prove necessity.
Let $u$ be a solution of the overdetermined problem
\begin{equation*}
\left\{
  \begin{array}{ll}
    \Delta_p^\mathcal{A} u=-1,~& \mathrm{in}~\Omega; \\
    u=0,~& \mathrm{on}~\partial\Omega;\\
    H_u(\nabla u)=c,~& \mathrm{on}~\partial\Omega.
  \end{array}
\right.
\end{equation*}
Set
$$H=H_u=\Big(\frac{\omega_n}{|L_{p,u}|}\Big)^\frac{1}{n}h_{\Gamma_pL_{p,u}}.$$
Then $H\in\mathcal{I}_p$ by Proposition \ref{33}. Since $\Delta_p^H=\Delta_p^\mathcal{A}$, the domain $\Omega$ and the function $u$ solve the overdetermined anisotropic problem
\begin{equation*}
\left\{
  \begin{array}{ll}
    \Delta_p^H u=-1,~& \mathrm{in}~\Omega; \\
    u=0,~& \mathrm{on}~\partial\Omega;\\
    H(\nabla u)=c,~& \mathrm{on}~\partial\Omega.
  \end{array}
\right.
\end{equation*}
Then, it follows from Theorem \ref{29} that, up to translations and dilations, $\overline{\Omega}$ is the Wulff shape of $\Big(\frac{\omega_n}{|L_{p,u}|}\Big)^\frac{1}{n}h_{\Gamma_pL_{p,u}}$, that is,
there are $z\in\rn$ and $c_1>0$ such that
$$\overline{\Omega}-z=c_1\bigcap_{y\in\rn}\Big\{x\in\rn:x\cdot y\leq \Big(\frac{\omega_n}{|L_{p,u}|}\Big)^\frac{1}{n}h_{\Gamma_pL_{p,u}}(y)\Big\}=c_1\Big(\frac{\omega_n}{|L_{p,u}|}\Big)^\frac{1}{n}\Gamma_pL_{p,u},$$
and
$$u(x)=\frac{(p-1)(r^{\frac{p}{p-1}}-H_0(x-z)^{\frac{p}{p-1}})}{pn^{\frac{1}{p-1}}}.$$
Set
$$g(s)=\frac{(p-1)r^\frac{p}{p-1}}{pn^\frac{1}{p-1}}(1-s^\frac{p}{p-1}),\qquad w(x)=u(x+z)=g\Big(\frac{H_0(x)}{r}\Big)=g\Big(\frac{1}{\rho_{\overline{\Omega}-z}(x)}\Big).$$
The function $g\in W^{1,p}([0,1])$ is strictly decreasing and satisfies $g(1)=0$, and $w\in W^{1,p}_0(\Omega-z)$. Moreover, $L_{p,u}=L_{p,w}$ by Lemma \ref{45}.
The set $K:=\overline{\Omega}-z$ is an origin-symmetric convex body. Applying Lemma \ref{34} to $w$ on $\Omega-z$ and using $\overline{\Omega-z}=K$, we obtain a constant $c_2>0$ such that
$$K=c_1\Big(\frac{\omega_n}{|L_{p,u}|}\Big)^\frac{1}{n}\Gamma_pL_{p,u}
=c_1\Big(\frac{\omega_n}{|L_{p,u}|}\Big)^\frac{1}{n}\Gamma_pL_{p,w}=c_1\Big(\frac{\omega_n}{|L_{p,u}|}\Big)^\frac{1}{n}\Gamma_p(c_2\Pi_p^\circ K).$$
Hence, there is a constant $c_3>0$ such that
$$K=c_3\Gamma_p\Pi_p^\circ K.$$
Theorem \ref{14}, applied to the convex body $K$, therefore implies that $K$ is an origin-symmetric ellipsoid, and hence $\Omega$ is an ellipsoid.

Since $K=\{H_0\leq r\}$ is an origin-symmetric ellipsoid, write $H_0(x)=|Ax|$. The normalization $H=H_u$ gives $\det A=1$, so $A=\psi\in\mathrm{SL}(n)$. Hence, with $z'=\psi z$,
$$u(x)=\frac{(p-1)(r^{\frac{p}{p-1}}-|\psi x-z'|^{\frac{p}{p-1}})}{pn^{\frac{1}{p-1}}},$$
for each $x\in\Omega$.
{\hfill$\Box$}
\section{Alexandrov's Theorem and Characterizations of Ellipsoids} \label{s4}
\subsection{The Affine Alexandrov Theorem and Heintze--Karcher Inequalities}
We begin with a regularity lemma that allows us to apply Theorem \ref{4} to the non-local affine mean curvature.
\begin{lemma}\label{lem:curvature-regularity}
Let $F$ be a set of finite perimeter with $C^2$ boundary and set $H=h_{\Gamma_1\Pi_1^\circ F}$. Then $H\in C^2(\rn\backslash\{0\})$, $H^2$ is uniformly convex, and
$$M_{\partial F}^H=\Big(\frac{|\Pi_1^\circ F|}{\omega_n}\Big)^{1/n}M_{\partial F}^\mathcal{A}.$$
\end{lemma}
\begin{proof}
By Lemma \ref{11}, $\Gamma_1\Pi_1^\circ F$ is of class $C^2_+$, and hence $H\in C^2(\rn\backslash\{0\})$. For $\xi\in S^{n-1}$, put $A_\xi=\nabla^2H(\xi)|_{\xi^\perp}$ and write $\nabla H(\xi)=H(\xi)\xi+z_\xi$, where $z_\xi\in\xi^\perp$. Since $H$ is positively homogeneous of degree one,
$$\nabla^2H(\xi)\xi=0,\qquad \nabla H(\xi)\cdot\xi=H(\xi)>0.$$
Since $A_\xi$ is positive definite on $\xi^\perp$, every $q=a\xi+w$ with $w\in\xi^\perp$ and $q\neq0$ satisfies
$$q^T\nabla^2(H^2)(\xi)q=2H(\xi)w^TA_\xi w+2\big(aH(\xi)+z_\xi\cdot w\big)^2>0,$$
because the first term is positive when $w\neq0$, while the second term is positive when $w=0$. By continuity and compactness, the smallest eigenvalue of $\nabla^2(H^2)$ has a positive lower bound on $S^{n-1}$. Positive homogeneity therefore implies that $H^2$ is uniformly convex. The curvature identity follows directly from the definitions of $M_{\partial F}^H$ and $M_{\partial F}^{\mathcal A}$.
\end{proof}

\noindent{\bf Proof of Theorem \ref{5}.} Let $\partial \Omega$ be a compact connected $C^2$ hypersurface without boundary embedded in Euclidean space and set $H=h_{\Gamma_1\Pi_1^\circ \Omega}$. By Lemma \ref{lem:curvature-regularity}, $H$ satisfies the hypotheses of Theorem \ref{4}, and $M_{\partial\Omega}^H$ is a constant positive multiple of $M_{\partial\Omega}^{\mathcal A}$. Hence, the constancy assumption in Theorem \ref{5} implies that $M_{\partial\Omega}^H$ is constant. Thus, by Theorem \ref{4}, $\overline{\Omega}$ is the Wulff shape of $H$ up to translations and dilations, that is, there are $z\in\rn$ and $c>0$ such that
$$c\Big(\overline{\Omega}+z\Big)=\bigcap_{y\in\rn}\Big\{x\in\rn:x\cdot y\leq h_{\Gamma_1\Pi_1^\circ \Omega}(y)\Big\}=\Gamma_1\Pi_1^\circ \Omega.$$
The translation invariance of the projection operator $\Pi$ (see Lemma \ref{41}) therefore implies that
$$\Gamma_1\Pi_1^\circ \Big(\overline{\Omega}+z\Big)=c\Big(\overline{\Omega}+z\Big).$$
Applying Theorem \ref{14} to the convex body $\overline{\Omega}+z$ shows that it is an origin-symmetric ellipsoid, and hence $\Omega$ is an ellipsoid.
{\hfill$\Box$}

The classical Heintze--Karcher inequality (see \cite{HK,ros}) states that
\begin{equation*}
  n|K|\leq \int_{\partial K}\frac{1}{M_{\partial K}}d\mathcal{H}^{n-1},
\end{equation*}
where $M_{\partial K}>0$ is the mean curvature of $\partial K$. Theorems \ref{4} and \ref{5} are essentially consequences of its anisotropic counterpart (Theorem \ref{16}), which also yields the affine Heintze--Karcher inequality in Theorem \ref{76}.

\begin{theorem}\label{16}
(\cite{JWXZ23}) Let $H$ be a norm on $\rn$ of class $C^2(\rn\backslash\{0\})$ such that $H^2$ is uniformly convex, and let $\partial \Omega$ be a compact connected $C^2$ hypersurface without boundary embedded in Euclidean space. If the anisotropic mean curvature $M_{\partial \Omega}^H>0$ on $\partial\Omega$, then
$$\int_{\partial\Omega}\frac{H\circ\nu_\Omega}{M_{\partial \Omega}^H}d\mathcal{H}^{n-1}\geq n|\Omega|,$$
with equality if and only if $\overline{\Omega}$ is the Wulff shape of $H$ up to translations and dilations.
\end{theorem}

As a direct consequence, we obtain the following affine Heintze--Karcher inequality:
\begin{theorem}\label{76}
Let $\partial \Omega$ be a compact connected $C^2$ hypersurface without boundary embedded in Euclidean space. If the non-local affine mean curvature $M_{\partial \Omega}^\mathcal{A}>0$ on $\partial\Omega$, then
$$\int_{\partial\Omega}\frac{(\omega_n/|\Pi_1^\circ\Omega|)^{1/n}h_{\Gamma_1\Pi_1^\circ\Omega}\circ\nu_\Omega}{M_{\partial \Omega}^\mathcal{A}}d\mathcal{H}^{n-1}\geq n|\Omega|,$$
with equality if and only if $\Omega$ is an ellipsoid.
\end{theorem}
\begin{proof}
With $H=h_{\Gamma_1\Pi_1^\circ\Omega}$, Lemma \ref{lem:curvature-regularity} shows that $H$ satisfies the hypotheses of Theorem \ref{16} and that $M_{\partial\Omega}^H$ is a constant positive multiple of $M_{\partial\Omega}^{\mathcal A}$. In particular, $M_{\partial\Omega}^H>0$, so the inequality follows from Theorem \ref{16}. If equality holds, Theorem \ref{16} gives $z\in\rn$ and $c>0$ such that $c(\overline{\Omega}+z)=\Gamma_1\Pi_1^\circ\Omega$. Translation invariance of $\Pi_1$ then yields
$$\Gamma_1\Pi_1^\circ(\overline{\Omega}+z)=c(\overline{\Omega}+z).$$
Applying Theorem \ref{14} to the convex body $\overline{\Omega}+z$ shows that it is an origin-symmetric ellipsoid, and hence $\Omega$ is an ellipsoid. Conversely, if $\Omega$ is an ellipsoid, choose $z_0$ such that $\overline{\Omega}+z_0$ is origin-symmetric. By Theorem \ref{14}, there is $c>0$ such that $\Gamma_1\Pi_1^\circ(\overline{\Omega}+z_0)=c(\overline{\Omega}+z_0)$. Translation invariance then gives $\Gamma_1\Pi_1^\circ\Omega=c(\overline{\Omega}+z_0)$. Hence $c(\overline{\Omega}+z_0)$ is precisely the Wulff shape of $H$, so equality follows from Theorem \ref{16}.
\end{proof}

\subsection{Mixed-Volume and Curvature Inequalities}
The Alexandrov--Fenchel inequality yields the following Minkowski-type inequalities for mixed volumes.
\begin{theorem}
(\cite{Sch14}) Let $K$ be a convex body with $C^2$ boundary and $L$ be a convex body with $C^2_+$ boundary. Then,
\begin{equation}\label{20}
  V_1(K,L)^2\geq|K|V_2(K,L)
\end{equation}
and
\begin{equation}\label{23}
  V_2(K,L)^{n-1}\geq V_1(K,L)^{n-2}|L|.
\end{equation}
Equality in \eqref{20} holds if and only if $K$ and $L$ are homothetic. For $n>2$, equality in \eqref{23} holds if and only if $K$ and $L$ are homothetic. When $n=2$, inequality \eqref{23} reduces to the identity $V_2(K,L)=|L|$, and therefore holds with equality for all $K$ and $L$.
\end{theorem}

\begin{lemma}\label{22}
(\cite{LYZ00}) Let $K$ be a convex body. Then,
$$V_1(K,\Gamma_1\Pi_1^\circ K)=\omega_n.$$
\end{lemma}

\noindent{\bf Proof of Theorem \ref{17}.} Let $K$ be a convex body with $C^2$ boundary. By Lemma \ref{11}, $\Gamma_1\Pi_1^\circ K$ is of class $C^2_+$. Thus, setting $L=\Gamma_1\Pi_1^\circ K$ in the Minkowski inequality \eqref{20}, one has
$$V_1(K,\Gamma_1\Pi_1^\circ K)^2\geq|K|V_2(K,\Gamma_1\Pi_1^\circ K),$$
that is,
$$\omega_n^2\geq|K|V_2(K,\Gamma_1\Pi_1^\circ K)$$
by Lemma \ref{22}.
By formula \eqref{49} and the definition of $M_{\partial K}^{\mathcal A}$,
$$V_2(K,\Gamma_1\Pi_1^\circ K)=\frac1n\Big(\frac{|\Pi_1^\circ K|}{\omega_n}\Big)^{1/n}\int_{\partial K}M_{\partial K}^{\mathcal{A}}h_{\Gamma_1\Pi_1^\circ K}(\nu_K(x))d\mathcal{H}^{n-1}(x).$$
Substituting this identity into the preceding inequality, we obtain
\begin{equation}\label{21}
  |K|\leq n\omega_n^{2+\frac1n}|\Pi_1^\circ K|^{-\frac1n}\bigg(\int_{\partial K}M_{\partial K}^{\mathcal{A}}h_{\Gamma_1\Pi_1^\circ K}(\nu_K(x))d\mathcal{H}^{n-1}(x)\bigg)^{-1}.
\end{equation}

If equality holds in \eqref{21}, then equality holds in \eqref{20}. Hence, by the equality condition in \eqref{20}, $K$ and $\Gamma_1\Pi_1^\circ K$ are homothetic. Thus, there exist $z\in\rn$ and $c>0$ such that $\Gamma_1\Pi_1^\circ K=cK+z$. Set $K_0=K+c^{-1}z$. By the translation invariance of the projection operator $\Pi$ (and hence of $\Pi_1$), it follows that
$$\Gamma_1\Pi_1^\circ K_0=\Gamma_1\Pi_1^\circ K=cK+z=cK_0.$$
Theorem \ref{14} therefore implies that $K_0$ is an origin-symmetric ellipsoid, and hence $K$ is an ellipsoid. Conversely, if $K$ is an ellipsoid, choose $z_0$ such that $K-z_0$ is origin-symmetric. By Theorem \ref{14}, there is $c>0$ such that $\Gamma_1\Pi_1^\circ(K-z_0)=c(K-z_0)$. Translation invariance gives $\Gamma_1\Pi_1^\circ K=cK-cz_0$, so $K$ and $\Gamma_1\Pi_1^\circ K$ are homothetic and equality holds in \eqref{20}, hence in \eqref{21}.
{\hfill$\Box$}

Finally, we show that the new affine curvature inequality is stronger than Petty's projection inequality.
We recall the Busemann--Petty centroid inequality as follows.
\begin{theorem}\label{24}
(\cite{Petty61}) Let $K$ be a convex body containing the origin in its interior. Then,
$$|\Gamma_1K|\geq|K|,$$
with equality if and only if $K$ is an origin-symmetric ellipsoid.
\end{theorem}

\begin{proposition}\label{25}
Let $K$ be a convex body with $C^2$ boundary. Then,
$$\bigg(\int_{\partial K}M_{\partial K}^{\mathcal{A}}h_{\Gamma_1\Pi_1^\circ K}(\nu_K(x))d\mathcal{H}^{n-1}(x)\bigg)^{n-1}\geq n^{n-1}\omega_n^{n-1-\frac1n}|\Pi_1^\circ K|^{\frac1n}.$$
For $n>2$, equality holds if and only if $K$ is an ellipsoid. For $n=2$, equality holds if and only if $K$ is a nondegenerate affine image of a planar convex body of constant width.
\end{proposition}
\begin{proof}
Set $P=\Pi_1^\circ K$ and $L=\Gamma_1P$. The Minkowski inequality \eqref{23}, Lemma \ref{22}, and Theorem \ref{24} yield
$$V_2(K,L)^{n-1}\geq V_1(K,L)^{n-2}|L|
=\omega_n^{n-2}|\Gamma_1P|\geq\omega_n^{n-2}|P|.$$
By formula \eqref{49} and the definition of $M_{\partial K}^{\mathcal A}$, this is equivalent to the asserted inequality.

If $n>2$, equality throughout implies that $K$ and $L$ are homothetic. Thus, there exist $z\in\rn$ and $c>0$ such that $\Gamma_1\Pi_1^\circ K=cK+z$. Setting $K_0=K+c^{-1}z$ gives $\Gamma_1\Pi_1^\circ K_0=cK_0$. Theorem \ref{14} therefore implies that $K_0$ is an origin-symmetric ellipsoid, and hence $K$ is an ellipsoid. Conversely, suppose that $K$ is an ellipsoid. Since projection bodies and polars of ellipsoids are ellipsoids, $P=\Pi_1^\circ K$ is an origin-symmetric ellipsoid. Theorem \ref{24} therefore gives $|\Gamma_1P|=|P|$. Choose $z_0$ such that $K-z_0$ is origin-symmetric. By Theorem \ref{14} and translation invariance, $L=\Gamma_1P=c(K-z_0)$ for some $c>0$. Hence $K$ and $L$ are homothetic, so equality holds in \eqref{23}. Thus both inequalities above are equalities.

If $n=2$, the first inequality is an identity because $V_2(K,L)=|L|$. Hence equality holds precisely when $|\Gamma_1P|=|P|$, that is, when $P$ is an origin-symmetric elliptic disk. Since polarity preserves centered elliptic disks, this is equivalent to $\Pi_1K=P^\circ$ being an elliptic disk. Let $J$ denote rotation through $\pi/2$. In the plane, the definition of $\Pi_1$ gives
$$\Pi_1K=\frac12J^{-1}(K-K).$$
Thus the equality condition is equivalent to $K-K$ being an elliptic disk. Writing $K-K=A(wB^2)$ with $A\in\mathrm{GL}(2)$ and $w>0$, this means precisely that $A^{-1}K$ has constant width $w$, namely that $K$ is a nondegenerate affine image of a planar convex body of constant width.
\end{proof}

Given a convex body $K$ with $C^2$ boundary, it follows from Theorem \ref{17} and Proposition \ref{25} that
$$n^{n-1}\omega_n^{n-1-\frac1n}|\Pi_1^\circ K|^{\frac1n}\leq \bigg(\int_{\partial K}M_{\partial K}^{\mathcal{A}}h_{\Gamma_1\Pi_1^\circ K}(\nu_K(x))d\mathcal{H}^{n-1}(x)\bigg)^{n-1}\leq (n\omega_n^{2+\frac1n})^{n-1}|\Pi_1^\circ K|^{-\frac{n-1}{n}}|K|^{1-n}.$$
These bounds imply Petty's projection inequality \eqref{18} and show that inequality \eqref{26} is stronger than Petty's inequality \eqref{18}.



\end{document}